\documentclass[a4paper,12pt]{amsart}

\usepackage{amsfonts}
\usepackage{amssymb}
\usepackage{ifthen}
\usepackage{graphicx}

\newtheorem{thm}{Theorem}
\newtheorem{cor}[equation]{Corollary}
\newtheorem{lem}{Lemma}
\newtheorem{prop}[equation]{Proposition}

\newtheorem{claim}{Claim}
\newtheorem{conj}[equation]{Conjecture}
\newtheorem{rem}{Remark}
\theoremstyle{definition}
\newtheorem{defn}{Definition}[section]
\newtheorem{example}{Example}[section]

\newtheorem{prob}[equation]{Problem}
\newtheorem{ques}[equation]{Question}

\newcounter {own}
\def\theown {\thesection       .\arabic{own}}

\newenvironment{pf}[1][]{%
 \vskip 2mm
 \noindent
 \ifthenelse{\equal{#1}{}}%
  {{\slshape Proof. }}%
  {{\slshape #1.} }%
 }%
{\qed\smallskip}

\newcounter{alphabet}
\newcounter{tmp}
\newenvironment{Thm}[1][]{\refstepcounter{alphabet}%
\bigskip%
\noindent%
{\bf Theorem \Alph{alphabet}}%
\ifthenelse{\equal{#1}{}}{}{ (#1)}%
{\bf .} \itshape}{\vskip 8pt}

\makeatletter
\newcommand{\Ref}[1]{\@ifundefined{r@#1}{}{\setcounter{tmp}{\ref{#1}}\Alph{tmp}}}
\makeatother
\newenvironment{Lem}[1][]{\refstepcounter{alphabet}%
\bigskip%
\noindent%
{\bf Lemma \Alph{alphabet}}%
{\bf .} \itshape}{\vskip 8pt}

\newcommand{\IC}{{\mathbb C}}
\newcommand{\ID}{{\mathbb D}}

\def\be{\begin{equation}}
\def\ee{\end{equation}}

\newcommand{\bee}{\begin{enumerate}}
\newcommand{\eee}{\end{enumerate}}

\newcommand{\blem}{\begin{lem}}
\newcommand{\elem}{\end{lem}}
\newcommand{\bthm}{\begin{thm}}
\newcommand{\ethm}{\end{thm}}
\newcommand{\bcor}{\begin{cor}}
\newcommand{\ecor}{\end{cor}}
\newcommand{\beg}{\begin{examp}}
\newcommand{\eeg}{\end{examp}}
\newcommand{\begs}{\begin{examples}}
\newcommand{\eegs}{\end{examples}}
\newcommand{\bdefe}{\begin{defn}}
\newcommand{\edefe}{\end{defn}}
\newcommand{\bprob}{\begin{prob}}
\newcommand{\eprob}{\end{prob}}
\newcommand{\bques}{\begin{ques}}
\newcommand{\eques}{\end{ques}}
\newcommand{\bei}{\begin{itemize}}
\newcommand{\eei}{\end{itemize}}

\newcommand{\bcon}{\begin{conj}}
\newcommand{\econ}{\end{conj}}
\newcommand{\bcons}{\begin{conjs}}
\newcommand{\econs}{\end{conjs}}
\newcommand{\bprop}{\begin{prop}}
\newcommand{\eprop}{\end{prop}}
\newcommand{\br}{\begin{rem}}
\newcommand{\er}{\end{rem}}
\newcommand{\brs}{\begin{rems}}
\newcommand{\ers}{\end{rems}}
\newcommand{\bo}{\begin{obser}}
\newcommand{\eo}{\end{obser}}
\newcommand{\bos}{\begin{obsers}}
\newcommand{\eos}{\end{obsers}}
\newcommand{\bpf}{\begin{pf}}
\newcommand{\epf}{\end{pf}}
\newcommand{\ba}{\begin{array}}
\newcommand{\ea}{\end{array}}
\newcommand{\beq}{\begin{eqnarray}}
\newcommand{\beqq}{\begin{eqnarray*}}
\newcommand{\eeq}{\end{eqnarray}}
\newcommand{\eeqq}{\end{eqnarray*}}

\newcommand{\ds}{\displaystyle}

\newcounter{minutes}
\divide\time by 60
\newcounter{hours}
\multiply\time by 60 \addtocounter{minutes}{-\time}
\begin{document}
%\begin{center}
%{ Communicated by W.R. Wogen}
%\end{center}
\bibliographystyle{amsplain}
\title [] {Area integral means, Hardy and weighted Bergman spaces of planar harmonic mappings}

%%%%%%%% BEGIN TIMESTAMP
\def\thefootnote{}
\footnotetext{ \texttt{\tiny File:~\jobname .tex,
          printed: \number\day-\number\month-\number\year,
          \thehours.\ifnum\theminutes<10{0}\fi\theminutes}
} \makeatletter\def\thefootnote{\@arabic\c@footnote}\makeatother
%%%%%%%% END TIMESTAMP

\author{SH. Chen}
\address{Sh. Chen, Department of Mathematics,
Hunan Normal University, Changsha, Hunan 410081, People's Republic
of China.} \email{shlchen1982@yahoo.com.cn}
\author{S. Ponnusamy}
\address{S. Ponnusamy, Department of Mathematics,
Indian Institute of Technology Madras, Chennai-600 036, India.}
\email{samy@iitm.ac.in}
\author{X. Wang${}^{~\mathbf{*}}$}
\address{X. Wang, Department of Mathematics,
Hunan Normal University, Changsha, Hunan 410081, People's Republic
of China.} \email{xtwang@hunnu.edu.cn}

\keywords{Harmonic mapping,  harmonic Hardy space, area integral means, and Bergmann spaces.\\
${}^{\mathbf{*}}$ Corresponding author}
%\date{\today  %June. 30, 09
%;  File: Ch-W-S13${}_{}$Har-Mean-Co.tex}

\begin{abstract}
In this paper, we investigate some properties of planar harmonic
mappings.  First, we generalize the main results in \cite{CPW3} and
\cite{HT}, and then discuss the relationship between area integral
means and harmonic Hardy spaces or harmonic weighted Bergman spaces.
At the end, coefficient estimates of mappings in weighted Bergman
spaces are obtained.
\end{abstract}

\maketitle
\pagestyle{myheadings}
\markboth{SH. Chen, S. Ponnusamy, X. Wang}{Area integral means, Hardy and weighted Bergman spaces}

\section{Introduction and preliminaries}\label{csw-sec1}

%$$\sup_{z,w\in\mathbb{B}^{n},z\neq w}\mathcal{L}_{f}(z,w)\leq2\pi\sqrt{n}\|f\|_{\mathcal{HB}}.
%$$
%Here$ 2\pi$ should be replaced by $\pi$. Am I right?

For each $r\in (0,1]$, we denote by $\ID_{r}$ the open disk $\{z\in
\IC:\, |z|<r\}$ and by $\ID$, the open unit disk $\ID_{1}$. The {\it
harmonic Hardy space} $\mathcal{H}_{h}^{p}(\mathbb{D})$ with
$0<p<\infty$ consists of all complex-valued functions $f$ harmonic
in  $\mathbb{D}$ (i.e. $f_{z\overline{z}}=0$ in $\ID$) for which
$$\|f\|_{p}:= \sup_{0<r<1}\big(I_{p}(r,f)\big)^{1/p}<\infty,
\quad
I_{p}(r,f)=\frac{1}{2\pi}\int_{0}^{2\pi}|f(re^{i\theta})|^{p}\,d\theta.
$$
The classical analytic Hardy space over the unit disk $\mathbb{D}$,
denoted usually by $\mathcal{H}^{p}(\mathbb{D})$, is obviously
contained in $\mathcal{H}_{h}^{p}(\mathbb{D}).$  We refer to
\cite{Clunie-Small-84, Du} for many basic analytic and geometric
properties of univalent harmonic mappings, in particular.

In this paper, we call a complex-value harmonic function as a
harmonic mapping. For a harmonic mapping $f$  in $\mathbb{D}$ and
$0\leq r<1$, the {\it generalized harmonic area function} $A_{h}(r)$
of $f$ is defined by (cf. \cite{CPW3})
%\be\label{eq1.1x}
$$A_{h}(r)=A_{h}(r,f)=\int_{\mathbb{D}_{r}}|\widehat{\nabla f}(z)|^{2}\,dA(z),
$$
where $dA$ denotes the normalized Lebesgue measure on $\mathbb{D}$,

$$\widehat{\nabla f}=(f_{z},f_{\overline{z}})~\mbox{and}~|\widehat{\nabla f}|=(|f_{z}|^{2}+|f_{\overline{z}}|^{2})^{1/2}.
$$
In particular, if $f$ is analytic in $\ID$, then we denote the
analytic  area function of $f$ by
$A(r):=A(r,f)=\int_{\mathbb{D}_{r}}|f'(z)|^{2}\,dA(z)$.

In \cite{HT,N}, the authors discussed the relationship between
(analytic) Hardy spaces and area functions. The  main  result in
\cite[Theorem 1]{HT} is as follows.

\begin{Thm}\label{Thm-1}
Let $f$ be analytic in $\mathbb{D}$. Then, if $1<p\leq 2$,
\be\label{eq3'}
f\in \mathcal{H}^{p}(\mathbb{D})\Rightarrow\int_{0}^{1}A^{p/2}
(r)\,dr<\infty,
\ee
while if $p>2$,
\be\label{eq4}
\int_{0}^{1}A^{p/2} (r)\,dr<\infty\Rightarrow f\in
\mathcal{H}^{p}(\mathbb{D}).
\ee
\end{Thm}
We refer to \cite{Du1,HL1,HL2,HT,MP1,MP2,MP3,MP4,N} for results
related to the theory of analytic Hardy spaces whereas for the
harmonic Hardy spaces, the readers may refer to \cite{CPW3,CPW4,D}.
In the context of recent investigation and interest on harmonic
mappings, it is natural to ask whether Theorem \Ref{Thm-1} continues
to hold in the setting of planar harmonic mappings over the unit
disk. In this note we show that the answer is yes.
%Theorem \Ref{Thm-1} continues to be true in the setting of planar harmonic
%mappings over the unit disk.

\begin{thm}\label{thm1}
Let $f$ be harmonic in $\mathbb{D}$. Then, if $1<p\leq 2$,
\be\label{eq3''} f\in
\mathcal{H}_{h}^{p}(\mathbb{D})\Rightarrow\int_{0}^{1}A_{h}^{p/2}
(r,f)\,dr<\infty, \ee while if $p>2$, \be\label{eq4'}
\int_{0}^{1}A_{h}^{p/2} (r,f)dr<\infty\Rightarrow f\in
\mathcal{H}_{h}^{p}(\mathbb{D}). \ee
\end{thm}

As an application of Theorem \ref{thm1}, we obtain the following
result.
%\Ref{Thm2}.

\begin{thm}\label{thm2}
Let $f\in\mathcal{H}_{h}^{p}(\mathbb{D})$. If $1<p\leq 2$, then $\ds
\lim_{r\rightarrow1-}(1-r)^{2/p}A_{h}(r,f)=0. $
\end{thm}

\br  {\rm Theorems \ref{thm1} and \ref{thm2} infer that the factor
$(1-r)^{\frac{\delta(2-p)}{2}}$ in \cite[Theorem 3]{CPW3} and the
one $(1-r)^{\frac{\delta(2-p)}{p}}$ in \cite[Theorem 4]{CPW3} are
redundant.
}\er

For a given real number $\alpha$, we consider the weighted area
measure $dA_{\alpha}^{\ast}(z)=(1-|z|^{2})^{\alpha}dA(z)$ on
$\mathbb{D}$, where $dA$ denotes the normalized area measure. For
$0<r<1$ and $0<p<\infty$, we define
%\be\label{eq10}
$$M_{p,\alpha}(r,f)=\left[\frac{1}{A_{\alpha}^{\ast}(\mathbb{D}_{r})}\int_{\mathbb{D}_{r}}
|f(z)|^{p}dA_{\alpha}^{\ast}(z)\right]^{1/p},
$$
where $f$ is harmonic in $\mathbb{D}$ and
$$A_{\alpha}^{\ast}(\mathbb{D}_{r})=\int_{\mathbb{D}_{r}}dA_{\alpha}^{\ast}(z).
$$

We call $M_{p,\alpha}(r,f)$, the {\it area integral means} of $f$ on $\mathbb{D}_r$.
%If $p=\infty$ and $f$ is harmonic in $\mathbb{D}$, then we define
%$$M_{\infty}(r,f)=\sup_{z\in\mathbb{D}_{r}}|f(z)|.
%$$

 It is well known that the measure
$A_{\alpha}^{\ast}$ is finite on $\mathbb{D}$ if and only if
$\alpha>-1.$ In the following, for $\alpha>-1$, we normalize the
measure $dA_{\alpha}^{\ast}$ by letting
$$dA_{\alpha}(z)=(1+\alpha)(1-|z|^{2})^{\alpha}dA(z).$$ For a harmonic mapping $f$  in $\mathbb{D}$, we denote
$$A_{\alpha}(\mathbb{D}_{r})=\int_{\mathbb{D}_{r}}dA_{\alpha}(z),
$$ where $\alpha>-1.$

For $\alpha>-1$ and $0<p\leq\infty$, the {\it weighted Bergman
space} $A_{h,\alpha}^{p}(\mathbb{D})$ consists of all harmonic
mappings $f$ on $\mathbb{D}$ such that
$$ \|f\|_{b^{p},\alpha}^{h}=
\begin{cases}
\displaystyle\left(\int_{\mathbb{D}}|f(z)|^{p}dA_{\alpha}(z)\right)^{1/p}<\infty
& \mbox{ if } p\in(0,\infty),\\
\displaystyle\sup_{z\in\mathbb{D}}|f(z)|<\infty &\mbox{ if }\,
p=\infty.
\end{cases}
$$

%In \cite{XZ}, the authors discussed some properties of the volume
%integral means on holomorphic mappings in $\mathbb{B}^{n}$.
Our next result
provides the relationship between area integral means
and harmonic Hardy spaces or harmonic weighted Bergman spaces.

\begin{thm}\label{thm3}
Suppose $1< p<\infty$, $\alpha$ is real, and $f$ is harmonic in
$\mathbb{D}$. Then, we have the following:
\begin{enumerate}
\item[{\rm \textbf{(a)}}] The function $M_{p,\alpha}(r,f)$ is
strictly increasing  in $[0,1)$ unless $f$ is constant.

\item[{\rm \textbf{(b)}}] For $\alpha>-1$,  $M_{p,\alpha}(r,f)$ is
bounded  in $[0,1)$ if and only if $f\in
A_{h,\alpha}^{p}(\mathbb{D})$.

\item[{\rm \textbf{(c)}}] For $\alpha\leq-1$,  $M_{p,\alpha}(r,f)$
is bounded in $[0,1)$ if and only if $f\in
\mathcal{H}_{h}^{p}(\mathbb{D}).$
\end{enumerate}
\end{thm}

Our final result concerns the coefficient estimate on
mappings in harmonic weighted Bergman spaces.

\begin{thm}\label{thm4}
For $1\leq p\leq\infty$, let $f\in A_{h,\alpha}^{p}(\mathbb{D})$ with
$$f(z)=\sum_{m=0}^{\infty}a_{m}z^{m}+\sum_{m=1}^{\infty}\overline{b}_{m}\overline{z}^{m}.$$
Then $|a_{0}|\leq \|f\|_{b^{p},\alpha}^{h}$, and for $m\geq 1$,
$$ |a_{m}|+|b_{m}|\leq\frac{4\|f\|_{b^{p},\alpha}^{h}}{\pi}\inf_{0<r<1}
\left \{\frac{1}{r^{m}\big[1-r^{\alpha+1}(2-r)^{\alpha+1}\big]^{1/p}}\right \}.
$$
In particular, if $\alpha=0$, then for $m\geq 1$,
$$|a_{m}|+|b_{m}|\leq\frac{4\|f\|_{b^{p},0}^{h}}{\pi}
\left (\frac{2}{pm}+1\right )^{m}\left(1+\frac{pm}{2}\right)^{\frac{2}{p}}.
$$
Moreover, if $\alpha=0$ and $p=\infty$, then
\be\label{eqthm4}
|a_{m}|+|b_{m}|\leq \frac{4\|f\|_{b^{\infty},0}^{h}}{\pi}.
\ee
The estimate \eqref{eqthm4} is sharp and the only extremal functions are
$$f_{m}(z)=\frac{2\gamma\|f\|_{b^{\infty},0}^{h} }{\pi}\arg \left ( \frac{1+\beta z^{m}}{1-\beta z^{m}}\right),
$$
where $ |\gamma|=|\beta |=1$, and $m\geq 1$.
\end{thm}

\section{Proofs of the main results}\label{csw-sec2}
We begin this section with the following two lemmas which are useful
in the proof of Theorem \ref{thm1}.

\begin{Lem}\label{Lem1}
Let $a, b\in[0,\infty)$ and $p\in[1,\infty)$. Then we have
$$a^{p}+b^{p}\leq(a+b)^{p}\leq 2^{p-1}(a^{p}+b^{p}).
$$
\end{Lem}

Lemma \Ref{Lem1} is well-known (see for instance \cite[Lemma 2.29 ]{S}).

\begin{lem}\label{lem1}
Let $f$ be a complex-valued continuously differentiable function
defined on $\mathbb{D}$ and $f=u+iv$, where $u$ and $v$ are
real-valued functions. Then for $z=x+iy\in\mathbb{D}$,
\be\label{eqs1}| f_{z}(z)|+|f_{\overline{z}}(z)|\leq |\nabla
u(x,y)|+|\nabla v(x,y)|, \ee where $\nabla u=(u_{x},u_{y})$ and
$\nabla v=(v_{x},v_{y})$.
\end{lem}
\bpf  A direct computation gives
$$f_{z}(z)=\frac{1}{2}(f_{x}(z)-if_{y}(z))~\mbox{and}~f_{\overline{z}}(z)=\frac{1}{2}(f_{x}(z)+if_{y}(z)),
$$
so that
$$f_{z}(z)=\frac{1}{2}[u_{x}(x,y)+v_{y}(x,y)+i(v_{x}(x,y)-u_{y}(x,y))]
$$
and
$$f_{\overline{z}}(z)=\frac{1}{2}[u_{x}(x,y)-v_{y}(x,y)+i(v_{x}(x,y)+u_{y}(x,y))].
$$
Applying the classical Cauchy-Schwarz inequality, we find that
\begin{eqnarray*}
|f_{z}(z)|&=&\frac{1}{2}\sqrt{\Big[\big(u_{x}(x,y)+v_{y}(x,y)\big)^{2}+\big(v_{x}(x,y)-u_{y}(x,y)\big)^{2}\Big]}\\
&\leq&\frac{1}{2}(|\nabla u(x,y)|+|\nabla v(x,y)|)
\end{eqnarray*}
and similarly,
\begin{eqnarray*}
|f_{\overline{z}}(z)|&=&\frac{1}{2}\sqrt{\Big[\big(u_{x}(x,y)-v_{y}(x,y)\big)^{2}+\big(v_{x}(x,y)+u_{y}(x,y)\big)^{2}\Big]}\\
&\leq&\frac{1}{2}(|\nabla u(x,y)|+|\nabla v(x,y)|),
\end{eqnarray*}
from which we easily obtain (\ref{eqs1}).
%$$| f_{z}(z)|+| f_{\overline{z}}(z)|\leq |\nabla u(x,y)|+|\nabla v(x,y)|.
%$$
%The proof of this Lemma is finished.
\epf

%We remark that the equality sign in (\ref{eqs1}) does not always
%hold. This can be seen from the following example.

\begin{example}
Let $f(z)=z^{2}+\overline{z}=u(x,y)+iv(x,y)$ in $\mathbb{D}$. Then
$u(x,y)=x^{2}+x-y^{2}$ and $v(x,y)=2xy-y$. It is easy to see that
$$|f_{z}(0)|+|f_{\overline{z}}(0)|=1 ~\mbox{ and }~ |\nabla u(0,0)|+|\nabla v(0,0)|=2.
$$
This observation shows that strict inequality in (\ref{eqs1}) is
possible.
\end{example}

\subsection*{Proof of Theorem \ref{thm1}}
We first prove the implication (\ref{eq3''}). Let $1<p\leq2$ and
$f=u+iv \in\mathcal{H}_{h}^{p}(\mathbb{D})$. Then $u$ and $v$ are
real harmonic functions in $\ID$. By Lemma \Ref{Lem1}, we deduce
that $u,\ v\in\mathcal{H}_{h}^{p}(\mathbb{D})$. Let $F_{1}$ and
$F_{2}$ be analytic functions defined on $\mathbb{D}$ such that
$\mbox{Re}\,F_{1}=u$ and $\mbox{Re}\,F_{2}=v$. Riesz' theorem (cf.
\cite[Theorem 4.1]{Du1}) shows that
$$\|F_{k}\|_{p}\leq\left(\frac{p}{p-1}\right)^{1/p}\|\mbox{Re}F_{k}\|_{p} ~\mbox{ for $k=1,2$}.
$$
which, in particular, implies that
$F_{k}\in\mathcal{H}^{p}(\mathbb{D})$ for $k=1,2$. By the
implication \eqref{eq3'} in Theorem \Ref{Thm-1}, it follows that
\be\label{eq6} \int_{0}^{1}A_{h}^{p/2}(r,F_{k})dr<\infty ~\mbox{ for
$k=1,2$}. \ee By calculations, we see that for $r\in (0,1)$,
\be\label{eq7}
A_{h}(r,F_{1})=\int_{\mathbb{D}_{r}}|F_{1}'(z)|^{2}dA(z)=\int_{\mathbb{D}_{r}}|\nabla
u(x,y)|^{2}dA(z)=A_{h}(r,u) \ee and similarly, \be\label{eq8}
A_{h}(r,F_{2})=\int_{\mathbb{D}_{r}}|F_{2}'(z)|^{2}dA(z)=\int_{\mathbb{D}_{r}}|\nabla
v(x,y)|^{2}dA(z)=A_{h}(r,v). \ee The inequalities (\ref{eq6}),
(\ref{eq7}), (\ref{eq8}) and Lemmas \Ref{Lem1} and \ref{lem1} yield
that
\begin{eqnarray*}
\int_{0}^{1}A_{h}^{p/2} (r,f)\,dr
&=&\int_{0}^{1}\left[\int_{\mathbb{D}_{r}}\big(|f_{z}(z)|^{2}+|f_{\overline{z}}(z)|^{2}\big)\,dA(z)\right]^{p/2}dr\\
&\leq&\int_{0}^{1}\left[\int_{\mathbb{D}_{r}}\big(|f_{z}(z)|+|f_{\overline{z}}(z)|\big)^{2}\,dA(z)\right]^{p/2}dr\\
&\leq&\int_{0}^{1}\left[\int_{\mathbb{D}_{r}}\big(|\nabla
u(x,y)|+|\nabla v(x,y)|\big)^{2}\,dA(z)\right]^{p/2}dr\\
&\leq&\int_{0}^{1}\left[\int_{\mathbb{D}_{r}}2\big(|\nabla
u(x,y)|^{2}+|\nabla v(x,y)|^{2}\big)dA(z)\right]^{p/2}dr\\
&\leq&2^{\frac{2p-1}{2}}\int_{0}^{1}\left[\left
(\int_{\mathbb{D}_{r}}|\nabla\
u(x,y)|^{2}dA(z)\right )^{p/2}\right. \\
&& \hspace{1.8cm} \left .+\left (\int_{\mathbb{D}_{r}}|\nabla
v(x,y)|^{2}\,dA(z)\right )^{p/2}\right]dr\\
&=&2^{\frac{2p-1}{2}}\int_{0}^{1}\left [A_{h}^{\frac{p}{2}}
(r,F_{1})+A_{h}^{\frac{p}{2}} (r,F_{2})\right ]dr\\
&<&\infty
\end{eqnarray*}
which proves the implication (\ref{eq3''}).

We next prove the implication (\ref{eq4'}). Let $p>2$ and $f$ be
harmonic in $\ID$. Then $f$  admits the canonical decomposition
$f=\phi +\overline{\psi}$, where $\phi $ and $\psi$ are analytic in
$\ID$ with $\psi(0)=0$. Then
$$\int_{\mathbb{D}_{r}}\big(|\phi'(z)|^{2}+|\psi '(z)|^{2}\big)\,dA(z)=A_{h} (r,f),
$$
which implies \be\label{eq9} \int_{0}^{1}A_{h}^{\frac{p}{2}}
(r,\phi)\,dr<\infty~\mbox{ and }~\int_{0}^{1}A_{h}^{\frac{p}{2}}
(r,\psi)\,dr<\infty. \ee By the implication \eqref{eq4} in Theorem
\Ref{Thm-1} and (\ref{eq9}), we conclude that $\phi,  \psi\in
\mathcal{H}^{p}(\mathbb{D})$. But then by the Minkowski inequality,
we deduce that
\begin{eqnarray*}
\big(I_{p}(r,f)\big)^{1/p}&\leq&\left(\frac{1}{2\pi}\int_{0}^{2\pi}\big(|\phi
(re^{i\theta})|+
|\psi (re^{i\theta})|\big)^{p}\,d\theta\right)^{1/p}\\
&\leq&\big(I_{p}(r,\phi )\big)^{1/p}+\big(I_{p}(r,\psi )\big)^{1/p},
\end{eqnarray*}
which yields that $\|f\|_{p}<\infty .$   %The proof of this theorem is finished.
\qed

\subsection*{Proof of Theorem \ref{thm2}}
It is not difficult to see that
$$(1-r)A_{h} ^{p/2}(r,f)\leq\int_{r}^{1}A_{h} ^{p/2}(\rho,f)\,d\rho, ~\mbox{ i.e. }~
%$$
%which implies
%$$
(1-r)^{2/p}A_{h} (r,f)\leq\left(\int_{r}^{1}A_{h} ^{p/2}(\rho,f)\,d\rho\right)^{2/p}.
$$
By the implication  (\ref{eq3''}) in Theorem \ref{thm1}, we conclude
$$\int_{0}^{1}A_{h}^{p/2} (r,f)\,dr<\infty
$$
from which we obtain that $\ds \lim_{r\rightarrow1-}(1-r)^{2/p}A_{h}(r,f)=0. $
%The proof of this theorem is finished.
\qed

\begin{lem}\label{lem-1}
Suppose that $f$ is harmonic on $\mathbb{D}$ and is constant in an open neighborhood of the origin.
%For $r\in(0,1)$, if
%$f$ in $\mathbb{D}_{r}$ is constant,
Then $f$ is constant throughout the unit disk $\mathbb{D}$.
\end{lem}
\bpf As every harmonic function $f$ in $\ID$ admits the representation
$$f(z)=a_{0}+\sum_{n=1}^{\infty}a_{n}z^{n}+\sum_{n=1}^{\infty}\overline{b}_{n}\overline{z}^{n},
%\quad  z\in \mathbb{D},
$$
we may assume that $f(z)= a_{0}$ in $\mathbb{D}_{r}$, for some $r\in (0,1)$. But then the
Parseval relation, for $0<\rho<r$, gives
$$|a_{0}|^{2}=\frac{1}{2\pi}\int_{0}^{2\pi}|f(\rho e^{i\theta})|^{2}d\theta=|a_{0}|^{2}+
\sum_{n=1}^{\infty}(|a_{n}|^{2}+|b_{n}|^{2})\rho^{2n},
$$
which obviously implies  $a_{n}=b_{n}=0$ for all $n\geq 1$. Thus, $f(z)\equiv a_{0}$ for
$z\in\mathbb{D}$.
\epf

Green's theorem (cf. \cite{CPW3,CPW4,P}) states that if $g\in
C^{2}(\mathbb{D})$, i.e., twice continuously differentiable on
$\ID$, then \be\label{eq1.2x}
\frac{1}{2\pi}\int_{0}^{2\pi}g(re^{i\theta})\,d\theta=g(0)+
\frac{1}{2}\int_{\mathbb{D}_{r}}\Delta g(z)\log\frac{r}{|z|}\,dA(z), \quad
\mbox{ $r\in (0, 1)$}.
\ee

\begin{lem}\label{lem-2}
Let $f$ be harmonic in $\mathbb{D}$. Then, for $p>1$, $I_{p}(r,f)$ is
a strictly increasing function of $r$ on $(0,1)$ unless $f$ is constant.
\end{lem}
\bpf By (\ref{eq1.2x}), we have
\begin{eqnarray*}
r\frac{d}{dr}I_{p}(r,f)&=&\frac{1}{2}\int_{\mathbb{D}_{r}}\Delta\big (|f(z)|^{p}\big )\,dA(z)\\
&=& p\int_{\mathbb{D}_{r}}\Big  [\big (\frac{p}{2}-1\big
)|f(z)|^{p-4} \big |f_{z}(z)\overline{f(z)}+f(z)\overline{f_{\overline{z}}(z)}\big |^{2}\\
&& \hspace{1.5cm} +|f(z)|^{p-2}|\widehat{\nabla f}(z)|^{2}\Big]\, dA(z)\\
&\geq&\frac{p(p-1)}{2}\int_{\mathbb{D}_{r}}\big(|f_{z}(z)|+|f_{\overline{z}}(z)|\big)^{2}|f(z)|^{p-2}\,
dA(z)\\
&\geq&0,
\end{eqnarray*}
which implies $I_{p}(r,f)$ is increasing on $r$ in $(0,1)$. Moreover, the last inequality
implies that $\frac{d}{dr}I_{p}(r,f)=0$ if and only if $f$ is
constant in $\mathbb{D}_{r}$. But then, in this case, Lemma \ref{lem-1} shows that $f$ is constant
on $\mathbb{D}$.
\epf

\subsection*{Proof of Theorem \ref{thm3}}
We first prove \textbf{(a)}. Since \be\label{eq11}
\int_{\mathbb{D}_{r}}|f(z)|^{p}dA_{\alpha}^{\ast}(z)=\int_{0}^{r}2\rho(1-\rho^{2})^{\alpha}I_{p}(\rho,f)d\rho,
\ee we see that \be\label{eq12}
\frac{d}{dr}\int_{\mathbb{D}_{r}}|f(z)|^{p}dA_{\alpha}^{\ast}(z)=2r(1-r^{2})^{\alpha}I_{p}(r,f).
\ee Simple calculations gives \be\label{eq13}
\frac{d}{dr}A_{\alpha}^{\ast}(\mathbb{D}_{r})=2r(1-r^{2})^{\alpha}.
\ee By (\ref{eq11}), (\ref{eq13}) and Lemma \ref{lem-2}, we have
$$I_p(r, f)-M_{p,\alpha}^{p}(r,f)=\frac{1}{A_{\alpha}^{\ast}(\mathbb{D}_{r})}
\int_{0}^{r}\left[\frac{d}{dt}I_p(t,
f)\right]A_{\alpha}^{\ast}(\mathbb{D}_{t})dt\geq0
$$
which implies \be\label{eq14} \int_{\mathbb{D}_{r}}
|f(z)|^{p}dA_{\alpha}^{\ast}(z)\leq
A_{\alpha}^{\ast}(\mathbb{D}_{r}) I_p(r, f). \ee By Lemma
\ref{lem-2}, we know that the equality holds in (\ref{eq14}) for
some $r$ only when $f$ is constant. By (\ref{eq12}), (\ref{eq14})
and computations, we conclude that
\begin{eqnarray*}
\frac{d}{dr}M_{p,\alpha}^{p}(r,f)&=&\frac{A_{\alpha}^{\ast}(\mathbb{D}_{r})\frac{d}{dr}\int_{\mathbb{D}_{r}}
|f(z)|^{p}dA_{\alpha}^{\ast}(z)-\int_{\mathbb{D}_{r}}
|f(z)|^{p}dA_{\alpha}^{\ast}(z)\frac{d}{dr}A_{\alpha}^{\ast}(\mathbb{D}_{r})}{A_{\alpha}^{\ast 2}(\mathbb{D}_{r})}\\
&=&\frac{2r(1-r^{2})^{\alpha}\Big[A_{\alpha}^{\ast}(\mathbb{D}_{r})I_p(r,
f)-\int_{\mathbb{D}_{r}}
|f(z)|^{p}dA_{\alpha}^{\ast}(z)\Big]}{A_{\alpha}^{\ast 2}(\mathbb{D}_{r})}\\
&\geq&0.
\end{eqnarray*}
Hence, the function $M_{p,\alpha}(r,f)$ is strictly increasing on
$r\in[0,1)$ unless $f$ is constant.

Next we prove \textbf{(b)}. We assume that $\alpha>-1$ and
$M_{p,\alpha}(r,f)$ is bounded. Then by \textbf{(a)}, we have
\be\label{eq16}\lim_{r\rightarrow1-}\left[\frac{1}{A_{\alpha}(\mathbb{D}_{r})}\int_{\mathbb{D}_{r}}
|f(z)|^{p}dA_{\alpha}(z)\right]=\int_{\mathbb{D}}|f(z)|^{p}dA_{\alpha}(z),\ee
which implies $f\in A_{h,\alpha}^{p}(\mathbb{D}).$ On the other
hand, if $f\in A_{h,\alpha}^{p}(\mathbb{D}),$ then the boundedness
of $M_{p,\alpha}(r,f)$ follows from (\ref{eq16}).

In order to prove \textbf{(c)}, we need some additional care.
\begin{claim}\label{cla1}
Suppose that $\alpha\leq-1$, $1\leq p<\infty$, and that $f$ is
harmonic in $\mathbb{D}$. Then %the integral
$$\int_{\mathbb{D}}|f(z)|^{p}dA_{\alpha}^{\ast}(z)
=\sup_{r\in(0,1)}\left\{\int_{\mathbb{D}_{r}}|f(z)|^{p}dA_{\alpha}^{\ast}(z)\right\}<\infty
~ \Longleftrightarrow ~f\equiv0
$$ %if and only if $f\equiv0$.
\end{claim}
\bpf Fix $\rho\in(0,1)$ so that $\rho<r<1$. Then \textbf{(a)} yields that
$$M_{p,\alpha}^{p}(R,f)\leq\frac{\int_{\mathbb{D}_{r}}
|f(z)|^{p}dA_{\alpha}^{\ast}(z)}{A_{\alpha}^{\ast}(\mathbb{D}_{r})}.
$$
It is not a difficult task to see that
$A_{\alpha}^{\ast}(\mathbb{D}_{r})\rightarrow \infty$ and
$$\int_{\mathbb{D}_{r}}|f(z)|^{p}dA_{\alpha}^{\ast}(z)\rightarrow \int_{\mathbb{D}}|f(z)|^{p}dA_{\alpha}^{\ast}(z)
~\mbox{ as $r\rightarrow 1- $},
$$
which gives  $M_{p,\alpha}(\rho,f)\equiv0$ for
each $\rho\in(0,1)$. Therefore, $f\equiv0$ and the proof of the claim is finished.
\epf

Finally, we prove \textbf{(c)}. We assume that $\alpha\leq-1$, $1<
p<\infty$ and $f$ is not identically zero. Then
$A_{\alpha}^{\ast}(\mathbb{D}_{r})\rightarrow \infty$ as
$r\rightarrow 1-$ and so Claim \ref{cla1} implies
$$\lim_{r\rightarrow1-}\int_{\mathbb{D}_{r}}|f(z)|^{p}dA_{\alpha}^{\ast}(z)=\int_{\mathbb{D}}|f(z)|^{p}dA_{\alpha}^{\ast}(z)=\infty.
$$
By \textbf{(a)} and calculations, we get
\begin{eqnarray*}
\sup_{0<r<1}M_{p,\alpha}^{p}(r,f)&=&\lim_{r\rightarrow1-}M_{p,\alpha}^{p}(r,f)\\
&=&\lim_{r\rightarrow1-}\frac{\int_{\mathbb{D}_{r}}|f(z)|^{p}(1-|z|^{2})^{\alpha}dA(z)}
{\int_{\mathbb{D}_{r}}(1-|z|^{2})^{\alpha}dA(z)}\\
&=&\lim_{r\rightarrow1-}\frac{2r(1-r^{2})^{\alpha}\frac{1}{2\pi}\int_{0}^{2\pi}|f(re^{i\theta})|d\theta}{2r(1-r^{2})^{\alpha}}\\
&=&\lim_{r\rightarrow1-}I_{p}(r,f)
%\\ &=&
=\|f\|_{p}^{p}
\end{eqnarray*}
and the proof of the theorem is complete.
\qed

\subsection*{Proof of Theorem \ref{thm4}} It is not difficult to show that for $p\in[1,\infty)$,
$|f|^{p}$ is subharmonic in $\mathbb{D}$. %Fix a point
%$z\in\mathbb{D}$ and let
%$$D=\{\zeta\in\mathbb{C}:\ |\zeta-z|<1-|z|\}.$$
Then for $z\in\mathbb{D}$ and $r\in[0,1-|z|)$, we have
$$|f(z)|^{p}\leq\frac{1}{2\pi}\int_{0}^{2\pi}|f(re^{i\theta}+z)|^{p}d\theta.$$
Integration gives
\begin{eqnarray*}
\big[1-|z|^{\alpha+1}(2-|z|)^{\alpha+1}\big]|f(z)|^{p}%&=&2\pi\int_{0}^{1-|z|}r|f(z)|^{p}dr\\
&\leq&\frac{1+\alpha}{\pi}\int_{0}^{2\pi}\int_{0}^{1-|z|}r(1-r^{2})^{\alpha}|f(z+re^{i\theta})|^{p}drd\theta\\
&\leq&\int_{\mathbb{D}}|f(z)|^{p}dA_{\alpha}(z)\, =\big(\|f\|_{b^{p},\alpha}^{h}\big)^{p}
\end{eqnarray*}
which implies
\be\label{eq1h}
|f(z)|\leq\frac{\|f\|_{b^{p},\alpha}^{h}}{\big[1-|z|^{\alpha+1}(2-|z|)^{\alpha+1}\big]^{\frac{1}{p}}}.
\ee
For $\zeta\in \mathbb{D}$ and $r\in(0,1)$, let $F(\zeta)=f(r\zeta)/r$. Then
$$F(\zeta)=\frac{a_{0}}{r} + \sum_{m=1}^{\infty}A_{m}\zeta^{m}+
\sum_{m=1}^{\infty}\overline{B}_{m}\overline{\zeta}^{m},
$$
where $A_m=a_{m}r^{m-1}~\mbox{ and }~ B_m=b_{m}r^{m-1}. $ Hence for
$\zeta\in\mathbb{D}$,
$$|F(\zeta)|\leq\frac{\|f\|_{b^{p},\alpha}^{h}}{r\big[1-r^{\alpha+1}(2-r)^{\alpha+1}\big]^{1/p}}=M(r).
$$

By (\ref{eq1h}), %it is not difficult to
we see that $|a_{0}|\leq \|f\|_{b^{p},\alpha}^{h}.$ It follows from
\cite[Lemma 1]{CPW2} that for $m\geq 1$,
$$|A_{m}|+|B_{m}|\leq\frac{4M(r)}{\pi}
$$ %~\mbox{ with $M(r)=\frac{\|f\|_{b^{p}}}{r(1-r)^{2/p}}.$}
which yields
$$ |a_{m}|+|b_{m}|\leq\frac{4\|f\|_{b^{p},\alpha}^{h}}{\pi}\inf_{0<r<1}
\left \{\frac{1}{r^{m}\big[1-r^{\alpha+1}(2-r)^{\alpha+1}\big]^{1/p}}\right \}.
$$
If $\alpha=0$, then
\begin{eqnarray*}
|a_{m}|+|b_{m}|&\leq&\frac{4\|f\|_{b^{p},0}^{h}}{\pi}\inf_{0<r<1}\left [\frac{1}{r^{m}(1-r)^{2/p}}\right ]
%\\
%&=&\frac{4}{\pi}\,\frac{1}{\ds \max_{ 0<r<1}\Big[r^{n}(1-r)^{2/p}\Big]}\\
%&=&
=\frac{4\|f\|_{b^{p},0}^{h}}{\pi}\left (\frac{2}{pm}+1\right )^{m}
\left(1+\frac{pm}{2}\right)^{\frac{2}{p}}.
\end{eqnarray*}
Since
$$\lim_{p\rightarrow\infty}\left (\frac{2}{pm}+1\right )^{m}\left(1+\frac{pm}{2}\right)^{\frac{2}{p}}=1,
$$
we conclude that \be\label{eq2.0} |a_{m}|+|b_{m}|\leq
\frac{4\|f\|_{b^{\infty},0}^{h}}{\pi}. \ee Thus, for $p=\infty$, the
estimate (\ref{eq2.0}) is sharp. By the subordination in the proof
of \cite[Theorem 1]{CPW-BMMSC2011},  we know that the only extreme
functions are
$$f_{m}(z)=\frac{2\gamma\|f\|_{b^{\infty},0}^{h}}{\pi}
\mbox{Im}\left(\log\frac{1+\beta z^{m}}{1-\beta z^{m}}\right) \quad
(|\gamma|=|\beta |=1),
$$
whose values are confined to $\ID_{\|f\|_{b^{\infty},0}^{h}} =
\{z\colon |z|<\|f\|_{b^{\infty},0}^{h}\}$.
%These functions lead to
%$$|a_{m}|+|b_{m}|=\frac{4\|f\|_{b^{\infty},0}^{h}}{\pi}.
%$$
%This completes the proof of our theorem.
\qed

\end{document}